%
%

\magnification 1200

\font\Bbb=msbm10
\def\BBB#1{\hbox{\Bbb#1}}

\def\HVir{{\cal L}}
\def\C{\BBB C}
\def\Z{\BBB Z}
\def\N{\BBB N}
\def\det{{\rm det}}
\def\char{{{\rm char}\hbox{\hskip 0.1cm}}}
\def\deg{{{\rm deg}\hbox{\hskip 0.1cm}}}

\def\hc{{h_I / c_{LI}}}
\def\o{{\bf 1}}
\def\d{\partial}
\def\Span{{\rm Span}}
\def\bw{{\overline w}}
\def\bu{{\overline u}}
\def\bv{{\overline v}}
\def\ide{{$I$-degree}}
\def\km{k}
\def\I{{\cal I}}

\def\ot{\otimes}

\def\LL{1.1}
\def\LI{1.2}
\def\II{1.3}
\def\phir{1.4}
\def\detf{1.5}

\def\He{2.1}
\def\B{2.2}
\def\dL{2.3}
\def\dI{2.4}
\def\BB{2.5}

\def\sing{1}
\def\ID{2}
\def\fund{3}
\def\vform{4}
\def\wV{5}
\def\MV{6}
\def\wfact{7}

\centerline
{\bf Representations of the twisted Heisenberg-Virasoro algebra at level 
zero.}

\

\centerline{
{\bf Yuly Billig}  
\footnote{*}{Research supported by the  Natural Sciences and
Engineering Research Council of Canada.}
\footnote{}{AMS subject classification: 17B68, 17B65.}
}

\

\

{\it Abstract.} We describe the structure of the irreducible highest weight modules
for the twisted Heisenberg-Virasoro Lie algebra at level zero. We prove that such
a module is either isomorphic to a Verma module or to a quotient of two Verma 
modules.

\

\

{\bf Introduction.}

\

 In this paper we study the structure of the irreducible representations for the twisted
Heisenberg-Virasoro Lie algebra $\HVir$ at level zero. This Lie algebra is the 
universal central extension of the Lie algebra of differential operators on a circle 
of order at most one:
$$ \left\{ f(t) {d \over dt} + g(t) \big| f,g \in \C[t,t^{-1}] \right\} .$$

The twisted Heisenberg-Virasoro algebra has an infinite-dimensional Heisenberg 
subalgebra and a Virasoro subalgebra. These subalgebras, however, do not form a 
semidirect product, but instead, the natural action of the Virasoro subalgebra on the 
Heisenberg subalgebra is twisted with a 2-cocycle (see (\LL)-(\II) for the precise
definition). 

The twisted Heisenberg-Virasoro algebra $\HVir$ has been 
studied by Arbarello et al. in [ACKP],
where a connection is established between the second cohomology of
certain moduli spaces of curves and the second cohomology of the Lie 
algebra of differential operators of order at most one.
Arbarello et al. also proved that when the central 
element of the Heisenberg subalgebra acts in a non-zero way, an irreducible highest 
weight module for $\HVir$ is isomorphic to the tensor product of an irreducible 
module for the Virasoro algebra and an irreducible module for the 
infinite-dimensional Heisenberg algebra.

The goal of the present paper is to study the case when the central element of the 
Heisenberg subalgebra acts trivially (level zero case). It turns out that the picture
in the level zero case is quite interesting and very different from the generic case
of non-zero level. Our main result (Theorem \sing \ below) states that either the 
Verma module itself is irreducible or the irreducible highest weight module is a 
quotient of two Verma modules. From this we immediately get the characters of the
irreducible modules for $\HVir$ at level zero.

Our work is motivated by the representation theory of the toroidal Lie algebras.
Level zero modules for the twisted Heisenberg-Virasoro algebra
appear as one of the ingredients in the construction of the vertex 
operator representations for the toroidal Lie algebras given in [B]. 

 We will denote by $\Z$ the set of integers and by $\N$ the set of natural 
numbers $\{ 1, 2, \ldots \}$.

{\bf Acknowledgements.} I am grateful to Victor Kac for bringing the paper [ACKP]
to my attention. I thank Yufeng Pei for pointing out a mistake in the original version of this paper.

\

{\bf 1. Twisted Heisenberg-Virasoro algebra.}

\

 We define the twisted Heisenberg-Virasoro algebra $\HVir$ as a Lie algebra with the 
basis
$$\big\{ L(n), I(n), C_L, C_{LI}, C_I, \big| n \in \Z \big\}$$
and the Lie bracket given by
$$ [ L(n) , L(m) ] = (n-m) L(n+m) + \delta_{n,-m} {n^3 - n \over 12} C_L 
,\eqno{(\LL)}$$
$$ [ L(n), I(m) ] = -m I(n+m) - \delta_{n,-m} (n^2 + n) C_{LI} ,
\eqno{(\LI)}$$
$$ [I(n), I(m) ] = n \delta_{n,-m} C_I, \eqno{(\II)}$$
$$[\HVir, C_L] = [\HVir, C_{LI}]  = [\HVir, C_I] = 0,  .$$

This Lie algebra has an infinite-dimensional Heisenberg subalgebra and a 
Virasoro subalgebra intertwined with the cocycle (\LI).
The twisted Heisenberg-Virasoro algebra $\HVir$ is the universal central extension of 
the Lie algebra $\{ f(t) {d \over dt} + g(t) | f,g \in\C[t,t^{-1}] \}$ of
differential operators of order at most one. The corresponding 
projection is given by $L(n) \mapsto -t^{n+1} {d \over dt}, $ \
$I(n) \mapsto t^n $.
The center of $\HVir$ is four-dimensional and is spanned by
$\{ I(0), C_L, C_{LI}, C_I \}$.

We are using the symbol $I$ because we may think of the infinite-dimensional 
Heisenberg algebra as the affinization of $gl_1(\C)$. In this 
interpretation $I$ is the identity matrix.   

 Introduce a $\Z$ grading on $\HVir$ by $\deg L(n) = \deg I(n) = n$ and
$\deg C_L = \deg C_{LI} = \deg C_I = 0$,
and decompose $\HVir$ with respect to this grading:
$\HVir = \HVir_- \oplus \HVir_0 \oplus \HVir_+$.

 Irreducible highest weight representations for $\HVir$ have been studied by
Arbarello et al. in [ACKP], however the case of irreducible representations
at level zero, i.e., when $C_I$ acts as zero, was not fully investigated 
in that paper.
It turns out that it is precisely this type of representations that
is needed for the construction of the modules for the toroidal Lie algebras
in [B]. The purpose of the present paper is to describe the 
structure of the irreducible 
modules for the twisted Heisenberg-Virasoro algebra at level zero. 
We are able to give a complete description (Theorem {\sing} below)
of these irreducible modules.

We begin by recalling the standard construction of the Verma modules.

Fix arbitrary complex numbers $h, h_I, c_L, c_{LI}, c_I$. Let $\C \o$ be a 
1-dimensional
$\HVir_0 \oplus \HVir_+$ module defined by $L(0) \o = h \o$, \
$I(0) \o = h_I \o$, \  $C_L \o = c_L \o$, \  $C_{LI} \o = c_{LI} \o$, \
$C_{I} \o = c_{I} \o$, \
$\HVir_+ \o = 0$.
As usual, the Verma module $M = M(h,h_I,c_L,c_{LI},c_I)$ is the induced 
module
$$ M(h,h_I,c_L,c_{LI},c_I) =
\hbox{\rm Ind}_{\HVir_0 \oplus\HVir_+}^{\HVir} (\C \o) \cong U(\HVir_-) 
\ot \o .$$

The module $M$ is $\Z$ graded by eigenvalues of the operator $L(0) - h$Id: 
\
$M = \mathop\oplus\limits_{n=0}^\infty M_n$ with $M_n = \{ v\in M | L(0) v 
= (n + h) v\}.$


 In order to understand the submodule structure of $M$, we need to study 
singular vectors
in $M$.
A non-zero homogeneous vector $v$ in a highest weight $\HVir$ module is 
called singular if $\HVir_+ v = 0$.

 Clearly, the highest weight vector $\o$ itself is singular, while every 
proper homogeneous submodule
of a highest weight module contains a singular vector which is not a 
multiple of the highest weight vector.

The key to the submodule structure of $M$ is the determinant formula
derived in [ACKP]. Let us briefly discuss this result.

The Lie algebra $\HVir$ has an anti-involution $\sigma$:
$$\sigma(L(n)) = L(-n), \quad \sigma(I(n)) = I(-n) - 2 \delta_{n,0}
C_{LI} ,$$
$$\sigma(C_L) = C_L, \quad \sigma(C_I) = C_I, \quad \sigma(C_{LI})
= -C_{LI} .$$

For a pair of modules
$$M^{(1)} = M(h,h_I,c_L,c_{LI},c_I), \quad
M^{(2)} = M(h, h_I - 2 c_{LI}, c_L, -c_{LI}, c_I),$$
this anti-involution induces a contragredient pairing
$$M^{(1)} \times M^{(2)} \rightarrow \C,$$
satisfying
$$ (x u | v) = (u | \sigma(x)v)$$ 
for all $x\in\HVir$, $u\in M^{(1)}$, $v\in M^{(2)}$ ,
and normalized by the condition $( \o | \o ) = 1$.
 
Clearly for $n\neq k$ we have $(M^{(1)}_n | M^{(2)}_k) = 0$.
Let us consider  the restriction of this pairing to 
$M^{(1)}_n \times M^{(2)}_n $.

If we fix a basis in the space $U_{-n}(\HVir_-)$ then we will get bases
in $M_n(h,h_I,c_L,c_{LI},c_I)$ simultaneously for all
$h,h_I,c_L,c_{LI},c_I \in\C$. Denote by $\det_n$ the determinant
of the contragredient pairing in this basis. 

If $\det_n = 0$ then both $M^{(1)}$ and $M^{(2)}$ have proper submodules with non-trivial $n$-th components. If $\det_n \neq 0$ 
for al $n$ then both $M^{(1)}$ and $M^{(2)}$ are irreducible.

Arbarello et al. established a formula
for $\det_n$ as a function of $h,h_I,c_L,c_{LI},c_I$
([ACKP], (6.7)).
In the case when $c_I = 0$, the determinant formula greatly
simplifies.
Here we present this reduction.
Define the numbers $p_2(n)$ by the generating series
$$ \sum\limits_{n=0}^\infty p_2(n) q^n =
\prod\limits_{k\geq 1} (1 - q^k)^{-2} ,$$
and let
$$\varphi_r = (h_I - (1+r)c_{LI})(h_I - (1-r)c_{LI}) .\eqno{(\phir)}$$
Then
$$ \det_n (h,h_I,c_L,c_{LI}) =
K_n \prod\limits_{1\leq s\leq r\leq n \atop 1\leq rs \leq n }
\varphi_{r,s}^{p_2 (n-rs)} , \eqno{(\detf)}$$
where
$$ \varphi_{r,s} = \left\{\matrix{
\varphi_r \varphi_s, \quad {\rm for} \quad r\neq s, \cr
\varphi_r , \quad {\rm for} \quad r=s, \cr}
\right. $$
and $K_n$ is a non-zero constant independent of $h,h_I,c_L,c_{LI}$
(in the notations of [ACKP] take $h_a = h_I - c_{LI}, c_3 = i c_{LI},
c_a = c_I = 0$).
Note that (\detf) is invariant under the substitution $h_I \mapsto h_I - 2 c_{LI}$, $c_{LI} \mapsto - c_{LI}$.

\

{\bf 2. Structure of the irreducible $\HVir$ modules at level zero.} 

\

The main result of the paper is the following

{\bf Theorem \sing.} 
{\it
Let $c_I = 0$ and $c_{LI} \neq 0.$

(a) If ${h_I \over c_{LI}} \notin \Z$ or  ${h_I \over c_{LI}} = 1$
then the $\HVir$ module $M = M(h,h_I,c_L,c_{LI},0)$ is irreducible.

(b) If ${h_I \over c_{LI}} \in \Z \backslash \{ 1 \}$ then
$M(h,h_I,c_L,c_{LI},0)$ possesses a singular vector $v\in M_p$,
where $p = \left| {h_I \over c_{LI}} - 1 \right|$.
The factor-module
$L = L(h,h_I,c_L,c_{LI},0) = M(h,h_I,c_L,c_{LI},0) / U(\HVir_-) v$ is 
irreducible
and its character is
$$ \char L = (1 - q^p) \prod\limits_{j\geq 1} (1 - q^j)^{-2} .$$
}

{\bf Proof.} 
Let us give a proof of part (a).
Clearly, if there exists a singular vector $v \in M_n$ for some $n > 0$
then $\det_n = 0$. Then the determinant formula (\detf) implies that 
$\varphi_m = 0$ for some $m\in\N$.
It follows from (\phir) that $\hc = 1 \pm m$. 
We conclude that if $\hc \notin \Z$ or $\hc = 1$ then
the Verma module $M$ does not possess a singular vector other than a vector
of the highest weight. Thus the Verma module $M$ is irreducible in this case.
This completes the proof of part (a) of the Theorem. 

To prove part (b), we will consider two cases: $1 - \hc \in \N$ and 
$\hc - 1 \in \N$. The proof in both cases is essentially the same, so we will
treat them in parallel.

Let us first outline the main idea of the proof. The Lie algebra $\HVir$ has
infinitely many Heisenberg subalgebras
$$ < L(n), \hbox{\hskip 0.2cm} 
I(-n),  \hbox{\hskip 0.2cm}
I(0) - (n+1) C_{LI} >_{n \neq 0} , 
\eqno{(\He)}$$
with the Lie bracket
$$ [ L(n), I(-n)] = n \left(  I(0) - (n+1) C_{LI} \right) .$$
The central element $ I(0) - (n+1) C_{LI}$ acts on $M =  M(h,h_I,c_L,c_{LI},0)$
in a non-zero way precisely when $\hc \neq 1 + n$. Thus in our case among the
Heisenberg subalgebras (\He), there will be one with
a degenerate action on $M$, and the rest will act non-degenerately.
We will exhibit a relation between the action of these Heisenberg subalgebras
and certain formal operations of taking partial derivatives in $M$ (see Lemma 
\fund \ below). The rest of the argument is reminiscent of the classical 
proof of irreducibility of a polynomial algebra as a module over a Heisenberg Lie 
algebra.

We will organize the proof of part (b) in a sequence of several lemmas.
For the rest of the paper we fix $p = 1 - \hc \in \N$ (resp. $p = \hc - 1 \in\N$).

From  (\phir) we get that $\varphi_r = - c_{LI}^2 (r-p)(r+p)$. Thus
$\varphi_p = 0$, while $\varphi_r \neq 0 $ for $r \neq p$, and it follows from
the determinant formula (\detf) that $\det_p = 0$, while $\det_{p-1} \neq 0$.
This implies the existence of a singular vector $v\in M_p$. Our goal is to
show that the submodule $V$ generated by this singular vector is the maximal 
submodule in $M$. To prove this, we need to study the properties of the singular
vector $v$ and of the submodule $V$.

Consider the following Poincar\'e -Birkhoff-Witt basis in $M = U(\HVir_-) \o$:
$$ \left\{  I(-m_1) \ldots I(-m_k) L(-n_1) \ldots L(-n_s) \o \right\}, \eqno{(\B)}$$
where $m_1 \geq \ldots \geq m_k > 0$, $n_1 \geq \ldots \geq n_s > 0$.
 
Note that the subalgebra $\HVir_-$ has one more $\Z$ grading by $I$-degree:
$$ \HVir_- = (\HVir_-)^I_0 \oplus (\HVir_-)^I_1, $$
where $I$-degree of $L(-n)$ is $0$, and  $I$-degree of $I(-n)$ is $1$.
We place a superscript $I$ in the notation of the graded component in order
to distinguish the grading by $I$-degree from the standard grading by the 
ordinary degree. This new grading on $\HVir_-$ induces a $\Z$ grading
on $U(\HVir_-)$ and also on the Verma module $M = U(\HVir_-) \o$:
$$ M = \mathop\oplus\limits_{j = 0}^\infty M^I_j .$$
The $I$-degree of a monomial in (\B) is $k$. 

 For a non-zero element $w\in M$ we will denote by $\bw$ its lowest 
non-zero homogeneous component with respect to \ide:
$$ w = \bw + \hbox{\rm terms of higher \ } I{\rm -degree.}$$

We define on $M$ the operations of formal partial derivatives 
${\d \over \d I(-m)}$, ${\d \over \d L(-n)}$. We set
$$ {\d I(-j) \over \d I(-m)} = \delta_{jm}, \quad
 {\d L(-j) \over \d I(-m)} = 0, \quad
 {\d \over \d I(-m)} \o = 0, $$
$$ {\d I(-j) \over \d L(-n)} = 0, \quad
 {\d L(-j) \over \d L(-n)} = \delta_{jn}, \quad
 {\d \over \d L(-n)} \o = 0, $$
and then define their action on monomials (\B) by the Leibnitz rule.
Finally, we extend these to $M$ by linearity. Clearly, these operations are
not canonical and depend on our choice of the basis.

\

{\bf Lemma \ID.} 
{\it
Let $w \in M_k^I$ and let $n > 0$. Then

(a) $I(n) w \in  M_k^I \oplus M_{k+1}^I$,

(b) $L(n) w \in M_{k-1}^I \oplus M_k^I$.
}

The proof of this Lemma is a simple application of the Poincar\'e -Birkhoff-Witt
argument and is left as an exercise.

 We will also need the following subspace in $M$:
$$ \I = \Span \left\{ I(-m_1) \ldots I(-m_k) \o \hbox{\hskip 0.2cm} 
\big| \hbox{\hskip 0.2cm}
m_1 \geq \ldots \geq m_k > 0 \right\} .$$

 The next Lemma exhibits the relation between the action of the Heisenberg
subalgebras (\He) and the formal partial derivatives.

{\bf Lemma \fund.} 
{\it
Let $w$ be a non-zero vector in $M$ expanded in the basis
(\B), and denote by $\km$ the \ide \ of its lowest component $\bw$.

(a) Suppose that $\bw \notin \I$. Let $n$ be the smallest integer such that 
$L(-n)$ occurs as a factor in one of the terms of $\bw$. Then the part of
$I(n) w $ of the \ide \  $\km$ is given by
$$ n (h_I + (n-1) c_{LI}) {d \bw \over \d L(-n)} .\eqno{(\dL)}$$

(b) Suppose that $\bw \in \I, \bw \notin \C\o$. Let $m$ be the maximal 
integer such that $I(-m)$ occurs as a factor in one of the terms of $\bw$.
Then the part of $L(m)w$ of the \ide \ $\km - 1$ is given by
$$ m (h_I - (m+1) c_{LI})  {d \bw \over \d I(-m)} \eqno{(\dI)}.$$
}
 
{\bf Proof.} Let us prove claim (a). By Lemma \ID (a), the part of $I(n) w$
of \ide \ $\km$ comes from $I(n) \bw$. Let
$$ x =  I(-m_1) \ldots I(-m_k) L(-n_1) \ldots L(-n_s) \o $$
be one of the monomials occurring in $\bw$. It is sufficient to establish
the claim of the Lemma for such a monomial. By our assumption, 
$n_1 \geq \ldots \geq n_s \geq n$. We have
$$ I(n) x = \sum_{i=1}^s I(-m_1) \ldots I(-m_k) L(-n_1) \ldots 
[I(n), L(-n_i)] \ldots L(-n_s) \o .$$
If $n_i > n$, we have $[I(n), L(-n_i)] = n I(-n_i + n)$, and the \ide 
\ of the corresponding term will be $\km + 1$, so such terms will not contribute
to the part of $I(n) x$ of \ide \ $\km$.

If $n_i = n$, then we have
$[I(n), L(-n)] = n(I(0) + (n-1) C_{LI})$, and so the contribution of these
terms is $ n (h_I + (n-1) c_{LI}) {\d x \over \d L(-n)} $. Combining these two
cases, we obtain claim (a).

 Let us now prove (b).  By Lemma \ID (b), the part of $L(m) w$
of \ide \ $\km - 1$ comes from $L(m) \bw$. 
Let
$$ y =  I(-m_1) \ldots I(-m_k) \o $$
be one of the monomials occurring in $\bw$. It is sufficient to establish  
the claim of the Lemma for such a monomial. By our assumption,
$m \geq m_1 \geq \ldots \geq m_\km$. We have
$$ L(m) y = \sum_{i=1}^\km I(-m_1) \ldots [L(m),I(-m_i)] \ldots I(-m_\km) \o .$$
If $m > m_i$ then  $[L(m), I(-m_i)] = m_i I(m-m_i)$, and the 
corresponding term vanishes since all $I$'s commute and $I(m-m_i) \o = 0$
when $m-m_i > 0$.

If $m = m_i$, then 
$[L(m), I(-m)] = m(I(0) - (m+1) C_{LI})$, and the contribution of these
terms will yield $ m (h_I - (m+1) c_{LI}) {\d y \over \d I(-m)} $.
The proof of the Lemma is now complete.

\

Note that the factor $ (h_I + (n-1) c_{LI})$ in (\dL) vanishes only when
$p = 1 - \hc \in\N$ and $n=p$. The factor $(h_I - (m+1) c_{LI})$ in (\dI)
is zero when $p = \hc -1 \in N$ and $m=p$. The partial derivatives
${\d \bw \over \d L(-n)}$ in (\dL) and ${\d \bw \over \d I(-m)}$ in (\dI)
are non-zero since by our assumptions $\bw$ involves $L(-n)$ in (a) and
$I(-m)$ in (b).   

In the following lemma we describe the decomposition of the singular vector 
$v\in M_p$ by \ide.

{\bf Lemma \vform.} 
{\it
The module $M$ possesses a singular vector $v \in M_p$
with $\bv = L(-p)\o$ (resp. $\bv = I(-p)\o$).
}

{\bf Proof.} We have already established the existence of a singular 
vector $v\in M_p$.
Denote the \ide \ of $\bv$ by $\km$.
Let us reason by contradiction and assume that $\bv$
is not a multiple of $L(-p) \o$ (resp. $I(-p) \o$). 
If $\bv \notin \I$, we apply Lemma \fund (a), and find $n \in \N$
such that $I(n) v \neq 0$ (note that in case  $p = 1 - \hc \in \N$ we have
$n \neq p$ due to our assumption that $\bv$ is not a multiple of $L(-p) \o$).
This contradicts to the fact that $v$ is a singular vector.    
If $\bv \in \I$ then we also get a contradiction in a similar way. We apply
Lemma \fund (b) to find $m \in \N$ such that $L(m) v \neq 0$ (note
that in case $p = \hc -1$, we have $m \neq p$ due to our assumption that
$\bv$ is not a multiple of $I(-p) \o$).    
Hence $\bv$ must be a multiple of $L(-p)\o$ (resp. $I(-p) \o$) and we  
rescale $v$ so that $\bv = L(-p)\o$ (resp. $\bv = I(-p)\o$).

\

{\bf Remark.} In fact it is possible to show that in case when $p = \hc -1\in\N$,
the singular vector $v$ belongs to $M_p \cap \I$.

\

{\bf Example.}
(i) If $\hc = 0$ then the singular vector
of degree 1 in $M(h,h_I,c_L,c_{LI},0)$ is $v = (L(-1) +
{h\over c_{LI}} I(-1) ) \o$. 

(ii) If $\hc = 2$ then the singular vector
of degree 1 in $M(h,h_I,c_L,c_{LI},0)$ is $v = I(-1)  \o$.

\

In order to prove Theorem \sing, we need to show that the submodule 
$V = U(\HVir_-)v$ generated by the singular vector $v \in M_p$ is the 
maximal submodule in $M$. To achieve this, we will need 
the following two corollaries to Lemma \vform.

\

{\bf Corollary \wV.} 
{\it
Let $w$ be a non-zero vector in the submodule 
$V = U(\HVir_-)v$, written in the basis (\B). Then there exist terms
in $\bw$, containing the factor $L(-p)$ (resp. $I(-p)$).
}

{\bf Proof.} Let $w = uv$, where $u \in U(\HVir_-)$. Since the universal 
enveloping algebra $U(\HVir_-)$ has no zero divisors, we have that 
$\bw = \bu$ $\bv$.
However by Lemma \vform, $\bv = L(-p)\o$ (resp. 
$\bv = I(-p)\o$). Thus $\bw = \bu L(-p) \o$ (resp.  $\bw = \bu I(-p) \o$).
Using the fact that the graded algebra $gr U(\HVir_-)$ associated with the
universal enveloping algebra $U(\HVir_-)$ is isomorphic to a polynomial
algebra, we conclude that all the terms in $\bw$ of the maximal length
(length of a monomial in (\B) is $s+k$) will contain a factor $L(-p)$ 
(resp. $I(-p)$). Thus we obtain the claim of the Corollary.

\

{\bf Corollary \MV.} 
{\it
The images of the vectors
$$ \left\{  I(-m_1) \ldots I(-m_k) L(-n_1) \ldots L(-n_s) \o \right\}, 
\eqno{(\BB)}$$
where $m_1 \geq \ldots \geq m_k > 0$, $n_1 \geq \ldots \geq n_s > 0$,
$n_i \neq p$ (resp. $m_i \neq p$), form the basis of the factor module $M/V$.
}

{\bf Proof.} By the previous Corollary, the vectors (\BB) are linearly independent 
modulo $V$. The character of the subspace in $M$ spanned by these vectors
coincides with the character of the factor module
$$\char M/V = (1 - q^p) \hbox{\hskip 0.2cm} \char U(\HVir_-) .$$
Thus the images of the vectors (\BB) under the projection $M \rightarrow M/V$
form the basis of $M/V$.

\

The next lemma is equivalent to the claim of part (b) of Theorem \sing. 
Our argument here will be quite similar to the one used in Lemma \vform.

\

{\bf Lemma \wfact.} 
{\it
Let $w \in M$. If $\HVir_+ w \subset V$ then 
$w \in \C\o \oplus V$.
}

{\bf Proof.} Without the loss of generality we may assume that $w$ is 
homogeneous. Also the statement of the Lemma will not change if we add to
$w$ a vector from $V$. Applying Corollary \MV, we may thus assume
that $w$ is a linear combination of vectors (\BB). 
We need to show that $w \in \C\o$.

If $\bw \notin\I$ then by Lemma 
\fund (a) there exists $n\in\N$ such that the part of $I(n)w$ of the 
lowest \ide \ is non-zero and belongs to the span of (\BB), because
the subspace spanned by vectors (\BB) is closed under the operations of
taking partial derivatives. However, by the assumption
of the Lemma, $I(n)w \in V$, which gives us a contradiction to Corollary \wV.

If $\bw \in\I, \bw\notin\C\o$, then by Lemma \fund (b) there exists
$m\in\N$ such that the part of the lowest \ide \ of $L(m)w$ is non-zero
and belongs to the span of (\BB). This again contradicts to Corollary \wV \ since
$L(m)w \in V$. This leaves us with the only possibility $\bw \in \C \o$ 
and hence $w\in\C\o$. The Lemma is proved.

\

Finally, to complete the proof of Theorem \sing (b), we note that by 
Lemma \wfact \ the only singular vectors in $L = M/V$ are 
multiples of the highest weight vector.
Thus $L$ is irreducible.

 The formula for the character follows from the obvious equalities
$$\char L = (1 - q^p) \hbox{\hskip 0.1cm} \char U(\HVir_-)$$
and
$$\char U(\HVir_-) = \prod\limits_{j\geq 1} (1 - q^j)^{-2} .$$

\

\

{\bf References:}

\

\noindent
[ACKP] E. Arbarello, C. De Concini, V.G. Kac, C. Procesi,
{\it Moduli spaces of curves and representation theory,}
Comm.Math.Phys., {\bf 117} (1988), 1-36.

\noindent
[B] Y. Billig, {\it Energy-momentum tensor for the toroidal Lie algebras,}
preprint.

\

\

\noindent
School of Mathematics and Statistics

\noindent
Carleton University

\noindent
1125 Colonel By Drive

\noindent
Ottawa, Ontario, K1S 5B6

\noindent
Canada

\noindent
e-mail: billig@math.carleton.ca

\end